\documentclass[10pt,twoside]{article}
\usepackage[latin1]{inputenc}
\usepackage{amsmath}
\usepackage{graphicx}
\usepackage{yhmath}
\usepackage[table]{xcolor}
\usepackage{mathrsfs} 
\usepackage{amssymb}
\usepackage{url}
\usepackage{makecell}
\usepackage{arydshln}
\usepackage{multirow}
\usepackage{arydshln}
\usepackage{mathtools}
\usepackage{stmaryrd}  

\input epsf
\def\limiten{\renewcommand{\arraystretch}{0.5}
\begin{array}[t]{c}\stackrel{}{\longrightarrow} \\
{\scriptstyle n\rightarrow
\infty}\end{array}\renewcommand{\arraystretch}{1}}

\def\limitepsn{\renewcommand{\arraystretch}{0.5}
\begin{array}[t]{c}\stackrel{a.s.}{\longrightarrow} \\
{\scriptstyle n \rightarrow
\infty}\end{array}\renewcommand{\arraystretch}{1}}

\def\limiteloin{\renewcommand{\arraystretch}{0.5}
\begin{array}[t]{c}\stackrel{{\cal D}}{\longrightarrow} \\
{\scriptstyle n\rightarrow
\infty}\end{array}\renewcommand{\arraystretch}{1}}

\def\limiteproban{\renewcommand{\arraystretch}{0.5}
\begin{array}[t]{c}\stackrel{{\cal P}}{\longrightarrow} \\
{\scriptstyle n\rightarrow
\infty}\end{array}\renewcommand{\arraystretch}{1}}

\def\equalpsn{\renewcommand{\arraystretch}{0.5}
\begin{array}[t]{c}\stackrel{a.s.}{=} \\
\end{array}\renewcommand{\arraystretch}{1}}

\numberwithin{equation}{section}

\newtheorem{thm}{Theorem}[section]

\newtheorem{lem}[thm]{Lemma}

\newtheorem{rmrk}[thm]{Remark}

\newcommand{\E}{\ensuremath{\mathbb{E}}}
\newcommand{\R}{\ensuremath{\mathbb{R}}}
\newcommand{\Z}{\ensuremath{\mathbb{Z}}}

\newcommand{\N}{\ensuremath{\mathbb{N}}}

\definecolor{grisclair}{gray}{0.9}

\renewcommand{\arraystretch}{.8}

\newcommand*\interior[1]{\overset{\mathsf{o}}{#1}}

\begin{document}
\title{\bf Some asymptotic results for time series model selection}
 \maketitle \vspace{-1.0cm}
  \begin{center}
    William KENGNE \footnote{Developed within the ANR BREAKRISK: ANR-17-CE26-0001-01 and the  CY Initiative of Excellence (grant "Investissements d'Avenir" ANR-16-IDEX-0008), Project "EcoDep" PSI-AAP2020-0000000013} 
 \end{center}

  \begin{center}
  { \it 
 THEMA, CY Cergy Paris Université, 33 Boulevard du Port, 95011 Cergy-Pontoise Cedex, France.\\ 
  william.kengne@cyu.fr  \\
  }
\end{center}

 \pagestyle{myheadings}
 \markboth{Some asymptotic results for time series model selection}{W. Kengne}

~~\\
\textbf{Abstract}:  We consider the model selection problem for a large class of time series models, including, multivariate count processes, causal processes with exogenous covariates. A procedure based on a general penalized contrast is proposed.
Some asymptotic results for weak and strong consistency are established.
The non consistency issue is addressed, and a class of penalty term, that does not ensure consistency is provided.
Examples of continuous valued and multivariate count autoregressive time series are considered.
\\ 


 \section{Introduction}\label{sec_intro}
Model selection is an important task in statistics and in many areas that processes data.
This issue, including time series model selection has gained increasing attention in the literature, given the large number of papers written in this direction in recent years.  We refer to the book of McQuarrie and Tsai (1998), the monograph of Rao and Wu (2001), the recent review paper of Ding {\it et al.} (2018)  and the references therein for an overview on this topic.
 
 \medskip

 Consider a trajectory $(Y_{1},\ldots,Y_{n})$, generated from a multivariate continuous valued or count time series  $Y=\{Y_{t},\,t\in \Z \}$, which belongs to a finite collection of models $\mathcal M$.
 In the sequel, it is assumed that $\mathcal M$ contains at least the true model $m^*$.
 For any model $m \in \mathcal M$, denote by $ \Theta(m)$ its parameter's space; the true model $m^*$ depends on $\theta^* \in \Theta(m^*)$. 
 Let $\Theta $ be a compact subset of $\R^d$ ($d \in \N$). We consider the following framework:
 \begin{itemize}
	\item each model $m \in \mathcal M$ is considered as a subset of $\{ 1,\ldots,d\}$ and denote by $|m|$ the dimension of  $m$ (i.e, $|m|=\#(m)$); 
		\item for  $m \in \mathcal M$, $\Theta(m)=\{  (\theta_i)_{1 \leq i \leq d} \in \Theta \text{ with } \theta_i=0 \text{ if } i \notin m \}$ is  the parameter's space of the model $m$; 
		\item $\mathcal{M}$ is considered as a subset of the power set of $\{ 1,\ldots,d\}$; i.e., $\mathcal{M} \subset \mathcal{P}(\{ 1,\ldots,d\})$.
\end{itemize}
 The aim is to select the "best model" (that we denote by  $\widehat m$) among the collection $\mathcal M$ such that it is "close" to $m^*$ for $n$ large enough.
 
\medskip

 For any model $m \in \mathcal M$, let $\widehat{\Phi} \big(Y_1,\cdots,Y_n;\theta\big)$ be a contrast function defined for all $\theta \in \Theta(m)$ by:
 \begin{equation}\label{C_def}
 \widehat{\Phi}_n \big(Y_1,\cdots,Y_n;\theta\big) = \sum_{t = 1}^n \widehat{\varphi}_t(\theta), 
  \end{equation} 
 where $\widehat{\varphi}_t(\cdot)$ depends on $Y_1,\ldots,Y_t$, such that, the minimum contrast estimator (MCE) given the model $m$,  is defined by 
 \begin{equation}\label{mce}
  \widehat{\theta}(m) \coloneqq  \underset{\theta\in \Theta(m)}{\text{argmin}} \Big( \widehat{\Phi}_n \big(Y_1,\cdots,Y_n;\theta\big) \Big).
  \end{equation}
In the sequel, the notation $\widehat{\Phi}_n (\theta  ) = \widehat{\Phi}_n \big(Y_1,\cdots,Y_n;\theta\big)$ is used. 
 For example, the contrast $\widehat{\Phi}$ could be the conditional likelihood, a quasi likelihood, a conditional least-squares or the density power divergence.
 Define the penalized criteria,
\begin{equation}\label{Cont_pen}
 \widehat{C}(m):= \widehat \Phi_n \big( (\widehat \theta(m) \big) + \kappa_n |m|, \text{ for all } m \in \mathcal{M} ,
 \end{equation}
 where $(\kappa_n)_{n \in \N}$ is a non negative sequence of the regularization parameter, satisfying $\kappa_n=o(n)$, and $|m|$ is the dimension of the model $m$. 
The selection of the "best" model $\widehat m$ among the collection $ \mathcal{M}$ is performed by minimizing the penalized contrast,
\begin{equation}\label{Estim_m}
  \widehat m :=  \underset{  m \in \mathcal{M}}{\text{argmin}} \left(\widehat{C}(m)\right).
  \end{equation} 

 \medskip
 
 The asymptotic properties usually considered for a model selection procedure are efficiency and consistency. 
 A procedure is said to be efficient according to a risk function, when its risk is asymptotically equivalent to that of the oracle. See for instance Shibata (1980), Karagrigoriou (1997), Ing and Wei (2005), Ing {\it et al.} (2012) for asymptotic efficiency results for model selection in $AR(\infty)$ type process. We also refer to the recent works of Hsu {\it et al.} (2019) and Bardet {\it et al.} (2021) for other results in a large class of linear and autoregressive processes.
 For consistency, it is assumed that the family of the competing model  contains at least the true $m^*$. A procedure is strongly (resp. weakly) consistent if the selected model $\widehat{m}$ converges $a.s.$ to $m^*$ (resp. the probability $P(\widehat{m} = m^*)$ approaches one) as the sample size tends to infinity. For example, see among others Hannan (1980), Tsay (1984)  for consistent results in ARMA-type model and the recent works of Bardet {\it et al.} (2020), Kengne (2021), Kamila (2021), Diop and Kengne (2022a, 2022b) for some results in a large class of autoregressive models, including models with exogenous covariates and integer-valued time series.
  
 \medskip
 
 In this new contribution, we study the consistent property for time series model selection. We do not set any model on the observations, thus, the theory developed here unifies many of the existing works and the results obtained apply to all models satisfying the conditions (\textbf{A1})-(\textbf{A4}) below. 
 Moreover, the contrast considered in this work is quite general and the existing procedures based on the likelihood, quasi likelihood, conditional least-squares or density power divergence can be seen as a specific case.   
The following issues are addressed.
\begin{itemize}
\item[(i)] Under the assumptions (\textbf{A1})-(\textbf{A4}) and some regular conditions on the parameter $(\kappa_n)$, the weak and the strong consistency of the procedure based on $\widehat{C}$ are established. In comparison with the existing results, besides the generality of our procedure, the conditions on the regularization parameter $\kappa_n$ are weaker (see also Remark \ref{rem_AC_X} and \ref{rem_MOD}).
\item[(ii)] We provide a class of $(\kappa_n)$ for which the weak and the strong consistencies fail. This issue generalizes the well known result of the non consistency of the AIC in time series model selection.
\item[(iii)] Applications to a large class of continuous valued and multivariate count autoregressive time series are considered.
\end{itemize}

  The rest of the paper is organized as follows. In Section 2, we set  some notations, assumptions and provide the asymptotic results of the proposed procedure. Application to a class of affine causal models with exogenous covariates is considered in Section 3. Section 4 focuses on a class of observation-driven multivariate count time series. Section 5 is devoted to the proofs of the main results.


  \section{Assumptions and asymptotic results}
  \subsection{Assumptions}
\noindent
 Throughout the sequel, we use the following norms:
{\em
\begin{itemize}
 \item $ \|x \| \coloneqq  \sum_{i=1}^{p} |x_i| $ for any $x \in \mathbb{R}^{p}$ (with $p \in \N$);
%
\item $ \|x \| \coloneqq \underset{1\leq j \leq q}{\max} \sum_{i=1}^{p} |x_{i,j}| $ for any matrix $x=(x_{i,j}) \in M_{p,q}(\R)$; where $M_{p,q}(\R)$ denotes the set of matrices of dimension $p\times q$ with coefficients in $\R$;
%
\item  $\left\|g\right\|_{\mathcal K} \coloneqq \sup_{\theta \in \mathcal K}\left(\left\|g(\theta)\right\|\right)$ for any  compact set $\mathcal K \subseteq  \Theta$ and 
function $g:\mathcal K \longrightarrow   M_{p,q}(\R)$;
%
 %
\item $\left\|Y\right\|_r \coloneqq \E\left(\left\|Y\right\|^r\right)^{1/r}$ for any random vector  $Y$ with finite $r-$order moments. 
\end{itemize}
}  
  
\medskip

  Consider a multivariate continuous or integer valued time series $Y=\{Y_{t},\,t\in \Z \}$, with the true model $m^*$, which depends on a parameter $\theta^* \in \Theta(m^*)$, and let $(Y_1,\cdots,Y_n)$ be a trajectory of $Y$. Denote by $\mathcal{F}_{t-1}=\sigma\left\{Y_{t-1},\ldots \right\}$ the $\sigma$-field generated by the whole past at time $t-1$.
 Let us set the following assumptions (see for instance Diop and Kengne (2021) for the change-point framework).
  \begin{enumerate} 
     \item [(\textbf{A1}):] The process $Y=\{Y_{t},\,t\in \Z \}$ is  stationary and ergodic.
    \item [(\textbf{A2}):] For any model $m \in \mathcal{M}$,  the function $\theta  \mapsto \widehat{\varphi}_t(\theta)$ (see (\ref{C_def})) is continuous on $\Theta(m)$; moreover, there exists a sequence of random function $(\varphi_t(\cdot))_{t \in \Z}$ such that, the mapping $\theta  \mapsto \varphi_t(\theta)$ is continuous on $\Theta(m)$ and for all $\theta \in \Theta(m)$, the sequence $( \varphi_t(\theta))_{t \in \Z}$ is  stationary and ergodic, satisfying:     
   \begin{equation}\label{eq_assump_app}
   \E \|\varphi_t(\theta) \|_{\Theta(m)} < \infty , ~
     \dfrac{1}{n}  \sum_{t=1}^n \| \widehat{\varphi}_t(\theta) 
     - \varphi_t(\theta)  \|_{\Theta(m)} 
     = o(1) ~ a.s.. 
  \end{equation}    
 In addition, there exists a unique $\theta^*(m)$ in $\Theta(m)$ such that the function $\theta  \mapsto \E[\varphi_0(\theta)]$ reaches its minimum in $\Theta(m)$ at $\theta^*(m)$ and $\theta^*$ is a unique point of $\Theta$ satisfying $\theta^* = \underset{\theta \in \Theta}{\text{argmin}}( \E[\varphi_0(\theta)]) $.
 In what follows, let us set $ \Phi_n (\theta) = \sum_{t = 1}^n \varphi_t(\theta)  $ for all $\theta \in \Theta$.
    \item [(\textbf{A3}):] For any $m \in \mathcal{M}$ and $ t =1,\cdots,n$, the function $\theta  \mapsto \widehat{\varphi}_t(\theta)$  is continuously differentiable on $\Theta(m)$; in addition, (\textbf{A2}) holds as well as 
  %
    the mapping $\theta  \mapsto \varphi_t(\theta)$ is continuously differentiable on $\Theta(m)$ and for all $\theta \in \Theta(m)$, the sequence $(\partial \varphi_t(\theta)/ \partial \theta)_{t \in \Z}$ is  stationary and ergodic, such as $ \E \| \partial \varphi_t(\theta) / \partial \theta \|^2_{\Theta(m)} < \infty$ and satisfying: 
   \begin{equation}\label{eq_assump_app_dif}
     \dfrac{1}{\sqrt{n}}  \sum_{t=1}^n \Big \| \dfrac{\partial}{\partial \theta} \widehat{\varphi}_t(\theta) 
     - \dfrac{\partial}{\partial \theta} \varphi_t(\theta) \Big \|_{\Theta(m)} 
     = o_P(1) ~ \text{ or } ~
  \end{equation} 
    \begin{equation}\label{eq_assump_app_dif_log2}
     \dfrac{1}{\sqrt{n \log \log n}}  \sum_{t=1}^n \Big \| \dfrac{\partial}{\partial \theta} \widehat{\varphi}_t(\theta) 
     - \dfrac{\partial}{\partial \theta} \varphi_t(\theta) \Big \|_{\Theta(m)} 
     = o(1) .
  \end{equation}  
  Furthermore, $\left( \big( \partial \varphi_t(\theta^*) /\partial \theta_i \big)_{i \in m},\mathcal{F}_{t}\right)_{t \in \mathbb{Z}}$ is a stationary ergodic, square integrable martingale difference sequence with covariance  $ G(m) = \Big( \E \Big[ \dfrac{\partial \varphi_0(\theta^*)}{\partial \theta_i}  \dfrac{\partial \varphi_0(\theta^*)}{\partial \theta_j}\Big] \Big)_{i,j \in m}$ assumed to exist and positive definite when $m^* \subseteq m$.  
 %
  %
  In the following, we set $G:=G(m^*)$. 
  \item [(\textbf{A4}):]  For any $m \in \mathcal{M}$ and $ t =1,\cdots,n$, the function $\theta  \mapsto \widehat{\varphi}_t(\theta)$  is two times continuously differentiable on $\Theta(m)$; moreover, (\textbf{A3}) holds as well as the function $\theta  \mapsto  \partial \varphi_t(\theta)/\partial \theta$ is continuously differentiable on $\Theta(m)$ and for all $\theta \in \Theta(m)$, 
the sequence $(\partial^2 \varphi_t(\theta)/ \partial \theta \partial \theta^T)_{t \in \Z}$ is  stationary and ergodic, satisfying:  
  %
   %
 \begin{equation}\label{eq_assump_app_dif2}
 \E \Big\|  \dfrac{\partial^2 \varphi_t(\theta)}{\partial \theta \partial \theta^T} \Big \|_{\Theta(m)} < \infty, ~ ~ 
    \Big \| \dfrac{1}{ n}  \sum_{t = 1}^n  \dfrac{\partial^2 \widehat{\varphi}_t(\theta)}{\partial \theta \partial \theta^T} - \E \Big(  \dfrac{\partial^2 \varphi_0(\theta)}{\partial \theta \partial \theta^T} \Big)  \Big \|_{\Theta(m)} 
     = o(1) ~ a.s.,
  \end{equation} 
  %
  %
  and the matrix  $F(m) =  \Big( \E \Big[ \dfrac{\partial^2 \varphi_0(\theta^*)}{\partial \theta_i \partial \theta_j}  \Big] \Big)_{i,j \in m}$ is assumed to exist and invertible when $m^* \subseteq m$. In the sequel, we set $F := F(m^*) $. 
  \end{enumerate}

 \begin{rmrk}\label{rm_assump}
 \begin{itemize}
 \item From the assumptions (\textbf{A1}) and (\textbf{A2}), it holds that, for any $m \in \mathcal{M}$, the MCE $ \widehat{\theta}(m)$ is consistent, that is 
 \begin{equation} \label{consist_normal_mce}
   \widehat{\theta}(m) \limitepsn \theta^*(m).
 \end{equation}
 Moreover, if $m^* \subset m$ (that is $\theta^* \in \Theta(m)$), then $\theta^*(m) = \theta^*$.
 \item Under (\textbf{A1})-(\textbf{A4}) and when $m^* \subseteq m$, standard arguments can be used to get the asymptotic normality of $\widehat{\theta}(m)$ (see also Lemma \ref{lem_a_n_theta_m_theta_star} below), that is 
  \begin{equation}\label{assymp_normal_mce}
  \sqrt{n}\left( \big(\widehat{\theta}(m) \big)_i- (\theta^*)_i\right)_{i \in m} \limiteloin \mathcal{N}_{|m|}(0,\Sigma)~ \text{ with } ~\Sigma(m) \coloneqq F(m)^{-1} G(m)  F(m)^{-1}. 
  \end{equation} 
  In the sequel, we set $\Sigma := F^{-1} G  F^{-1}$.
 \end{itemize}
 \end{rmrk} 
   
  \medskip
 
  \noindent
The examples of Section \ref{Sec_class_AC} and \ref{Sec_class_OD} show that, the assumptions (\textbf{A1})-(\textbf{A4}) hold for many classical models.

 \subsection{Asymptotic results} 
 The following theorem provides sufficient conditions for weak  consistency of the model selection procedure (\ref{Estim_m}), as well as some conditions for which the consistency fails.
\begin{thm}\label{th1_mce}
Assume (\textbf{A1}), (\textbf{A2}), (\textbf{A3}) with (\ref{eq_assump_app_dif}), (\textbf{A4}) and $\theta^* \in \overset{\circ}{\Theta}$. The following hold:
\begin{enumerate}
\item If  $k_n \limiten \infty$, then, $\widehat{m} \limiteproban m^*$.
\item If $ \displaystyle{\limsup_{n \rightarrow \infty} \kappa_n < \infty}$, then, the weak consistency of $\widehat{m}$ fails, unless $m^* \not\subset m$ for all $m \in \mathcal{M} \setminus \{m^*\}$.
\end{enumerate} 
\end{thm}
It is well known that the AIC asymptotically tends to choose overfitting models and is not consistent (see for instance Shibata (1976)); the second part of Theorem \ref{th1_mce} extends such result.     
 The next theorem provides sufficient conditions for strong consistency of the procedure (\ref{Estim_m}), as well as some conditions for which the strong consistency fails. 
  \medskip

  \begin{thm}\label{th2_mce}
Assume (\textbf{A1}), (\textbf{A2}), (\textbf{A3}) with (\ref{eq_assump_app_dif_log2}), (\textbf{A4}) and $\theta^* \in \overset{\circ}{\Theta}$. We have the following:
\begin{enumerate}
\item There exists a constant $c_1>0$ such that if $ ~\underset{n\rightarrow \infty}{\liminf}  (\kappa_n/\log \log n)>c_1$, then 
$
 \widehat m  \limitepsn m^* . 
 $ 
\item Assume $\exists \vartheta >0$ such that $G(m) = \vartheta F(m)$ for any $m \in \mathcal{M}$ satisfying $m^* \subseteq m$.  There exists a constant $c_2>0$ such that if $ ~\underset{n\rightarrow \infty}{\limsup}  (\kappa_n/\log \log n) < c_2$, then the strong consistency of $\widehat{m}$ fails, unless $m^* \not\subset m$ for all $m \in \mathcal{M} \setminus \{m^*\}$.
\end{enumerate} 
\end{thm}
Let us point out that, the condition $G(m) = \vartheta F(m)$ holds for many classical models. For $AR(\infty)$ (including $ARMA$) processes, it holds with $\vartheta=2$ and for $ARCH(\infty)$ (including $GARCH$) processes, it holds with $\vartheta= \E \xi_0^4 -1$ (when $\E \xi_0^4 < 1$) where $(\xi_t)$ is the innovations process of the model (see for instance Bardet and Wintenberger (2009)).


 \section{Application to a class of affine causal models with exogenous covariates}\label{Sec_class_AC}
 Let $X_t=(X_{1,t},X_{2,t},\ldots,X_{d_{x}, t}) \in \R^{d_{x}}$ a vector of covariates, with $d_{x} \in \N$.  Consider the class of affine causal models with exogenous covariates (see Diop and Kengne (2022a)) given by, 
 
 \medskip
 
 \noindent \textbf{Class} $\mathcal {AC}$-$X(M_\theta,f_\theta):$ A process $\{Y_{t},\,t\in \Z \}$ belongs to $\mathcal {AC}$-$X(M_\theta,f_\theta)$ if it satisfies:
   \begin{equation}\label{Model_AC_X} 
     Y_t =M_\theta(Y_{t-1}, Y_{t-2}, \ldots; X_{t-1},X_{t-2},\ldots)\xi_t + 
      f_\theta(Y_{t-1}, Y_{t-2}, \ldots; X_{t-1},X_{t-2},\ldots),
   \end{equation}
  where 
   $M_\theta,~f_\theta :  \R^{\N} \times (\R^{d_x})^{\N} \rightarrow \R$ are two measurable functions 
   and assumed to be known up to the parameter $\theta$, which belongs to a compact subset $ \Theta\subset \R^d$ ($d \in \N$); 
   and
   $(\xi_t)_{t \in \Z}$ is a sequence of zero-mean independent, identically distributed (\textit{i.i.d}) random variable satisfying $\E(\xi^r_0) < \infty$ for some $r \geq 2$ and $\E(\xi^2_0) =1$. 
   Note that, if the covariates is absent (that is $X_t \equiv C$ for some constant $C$), then, (\ref{Model_AC_X}) reduces to the classical affine causal models that has been studied among other by  Bardet and Wintenberger (2009), Bardet {\it et al.} (2012), Bardet {\it et al.} (2020)). One can see that, the  ARMAX, TARX, GARCH-X, APARCH-X (see Francq and Thieu (2019)), models belong to the class $\mathcal {AC}$-$X(M_\theta,f_\theta)$.
  Diop and Kengne (2022a) introduced other examples such as APARCH-X$(\delta,\infty)$, ARX($\infty$)-ARCH($\infty$), which belong to $\mathcal {AC}$-$X(M_\theta,f_\theta)$.  
 
    \medskip
       
    \noindent
  Consider a trajectory $(Y_1,\ldots,Y_n)$ generated from a process $(Y_t)_{t \in \Z}$ that belongs to $\mathcal {AC}$-$X(M_{\theta^*},f_{\theta^*})$, where $\theta^* \in \Theta \subset \R^d$ is the true parameter.
  Let $\mathcal M$ be a finite collection of the competing models that contains at least the true model $m^*$ corresponding to the parameter $\theta^*$. The aim is to select the "best model" among the collection $\mathcal M$. To this end,   we would like to carry out the general procedure developed in Section \ref{sec_intro}.
  
   \medskip
   
For any segment $\theta \in \Theta$, define the contrast function based on the conditional Gaussian quasi (log)likelihood given (up to an additional constant) by
%
%

%
 \begin{equation}\label{def_Phi.theta_AC_X} 
 \widehat  \Phi_n(\theta)=  \frac{1}{2}\sum_{t \in T} \widehat \varphi_t(\theta)
 ~~~ \text{with}~~~
        \widehat \varphi_t(\theta) =  \dfrac{(Y_t-\widehat f^t_\theta)^2}{\widehat H^t_\theta}+\log(\widehat H^t_\theta),
\end{equation}      
where $\widehat f^t_\theta := f_\theta(Y_{t-1},\ldots,Y_{1},0,\ldots ; X_{t-1},\ldots,X_{1},0,\ldots)$, 
$\widehat M^t_\theta := M_\theta(Y_{t-1},\ldots,Y_{1},0,\ldots; X_{t-1},\ldots,X_{1},0,\ldots)$ and $\widehat H^t_\theta := \big(\widehat M^t_\theta \big)^2$.
 %
 %
 Therefore, the MCE,  given a model $m$ is defined by
 \begin{equation}\label{mce_AC_X} 
 \widehat{\theta}(m)= \underset{\theta\in \Theta(m)}{\text{argmin}} \big(\widehat \Phi_n(\theta) \big).
\end{equation}

  \medskip
 
 Diop and Kengne (2022a) imposed the following Lipschitz-type conditions on the function $f_\theta$, $M_\theta$ or $M_\theta^2$ in order to study the stability properties of the class $\mathcal {AC}$-$X(M_\theta,f_\theta)$.
  In the sequel, denote by $0$ the null vector of any vector space. 
 For $\Psi_\theta = \theta$ or $M_\theta$  any compact set $\mathcal K \subseteq \Theta$, consider the assumption
 
   \medskip
   
    \noindent 
    \textbf{Assumption} \textbf{A}$_i (\Psi_\theta,\mathcal K)$ ($i=0,1,2$):
    For any $(y,x) \in \R^{\infty}\times (\R^{d_x})^{\infty}$, the function $\theta \mapsto \Psi_\theta(y)$ is $i$ times continuously differentiable on $\mathcal K$  with $ \big\| \frac{\partial^i \Psi_\theta(0)}{\partial \theta^i}\big\|_{\mathcal K}<\infty $; 
    and
      there exists two sequences of non-negative real numbers $(\alpha^{(i)}_{k,Y}(\Psi_\theta,\mathcal K))_{k\geq 1} $ and $(\alpha^{(i)}_{k,X}(\Psi_\theta,\mathcal K))_{k\geq 1}$ satisfying:
     $ \sum\limits_{k=1}^{\infty} \alpha^{(i)}_{k,Y}(\Psi_\theta,\mathcal K) <\infty$, $ \sum\limits_{k=1}^{\infty} \alpha^{(i)}_{k,X}(\Psi_\theta,\mathcal K) <\infty$  for $i=0, 1, 2$;
   such that for any  $(y,x), (y',x') \in \R^{\infty} \times (\R^{d_x})^{\infty}$,
  \[  \Big \| \frac{\partial^i \Psi_\theta(y,x)}{ \partial \theta^i}-\frac{\partial^i \Psi_\theta(y',x')}{\partial\theta^i} \Big \|_{\mathcal K}
  \leq  \sum\limits_{k=1}^{\infty}\alpha^{(i)}_{k,Y}(\Psi_\theta,\mathcal K) |y_k-y'_k|+\sum\limits_{k=1}^{\infty}\alpha^{(i)}_{k,X}(\Psi_\theta,\mathcal K)\|x_k-x'_k\|, 
  \]
where $\| \cdot\|$ denotes a vector, matrix norm.
\medskip
   
    \noindent
The next assumption is set on the function $H_\theta=M^2_\theta$ in the cases of ARCH-X type process.
\medskip
   
    \noindent
 \textbf{Assumption} \textbf{A}$_i (H_\theta,\mathcal K)$ ($i=0,1,2$): Assume $f_\theta=0$ and that 
  for any $(y,x) \in \R^{\infty}\times (\R^{d_x})^{\infty}$, the function $\theta \mapsto H_\theta=M^2_\theta(y)$ is $i$ times continuously differentiable on $\mathcal K$  with $ \big\| \frac{\partial^i H_\theta(0)}{\partial \theta^i}\big\|_{\mathcal K}<\infty $; and there exists two sequences of non-negative real numbers $(\alpha^{(i)}_{k,Y}(H_\theta,\mathcal K))_{k\geq 1} $ and $(\alpha^{(i)}_{k,X}(H_\theta,\mathcal K))_{k\geq 1} $ satisfying:
     $ \sum\limits_{k=1}^{\infty} \alpha^{(i)}_{k,Y}(H_\theta,\mathcal K) <\infty $, $ \sum\limits_{k=1}^{\infty} \alpha^{(i)}_{k,X}(H_\theta,\mathcal K) <\infty $  for $i=0, 1, 2$;
   such that for any  $(y,x), (y',x') \in \R^{\infty} \times (\R^{d_x})^{\infty}$,
  \[  \Big \| \frac{\partial^i H_\theta(y,x)}{ \partial \theta^i}-\frac{\partial^i H_\theta(y',x')}{\partial\theta^i} \Big \|_{\mathcal K}
  \leq  \sum\limits_{k=1}^{\infty}\alpha^{(i)}_{k,Y}(H_\theta,\mathcal K) |y^2_k-{y'}^2_k|+\sum\limits_{k=1}^{\infty}\alpha^{(i)}_{k,X}(H_\theta,\mathcal K) \|x_k-x'_k\| . 
  \]
  
\medskip

  \noindent
Also, we impose an autoregressive-type structure on the covariates:
 \begin{equation}\label{Process_X_AC_X}
 X_t=g(X_{t-1},X_{t-2},\ldots;\eta_t),
\end{equation}
 where $(\eta_t)_{t \in \Z}$ is a sequence of centered random variables such as $(\eta_t, \xi_t)_{t \in \Z}$  is {\it  i.i.d}
  and $g$ is a $\R^{d_x}$-valued function such that
\begin{equation} \label{exp_Lip_cov_AC_X}  
   \E\left[\left\|g(0, \eta_0) \right\|^r\right]<\infty
   ~ \text{ and }  ~ 
   \left\|g(x; \eta_0)-g(x'; \eta_0) \right\|_r \leq \sum\limits_{k=1}^{\infty} \alpha_k(g) \left\| x_k-x'_k \right\|
   ~ \text{ for all } x, x'  \in (\R^{d_{x}})^{\infty},
  \end{equation}
   for some $r \geq 1$ and non-negative sequence $(\alpha_k(g))_{k \geq 1}$  fulfilling $\sum\limits_{k=1}^{\infty} \alpha_k(g) <1 $.
   
   \medskip
   
    \noindent
    In this section, we assume that (\ref{Process_X_AC_X}) and (\ref{exp_Lip_cov_AC_X}) hold for some $r\geq 1$.
  Define the set 
\begin{multline*}\label{Set_Theta(r)_AC_X}
\Theta(r) = \Big\{
 \theta \in \R^d \, \big / \, \textbf{A}_0 (f_\theta,\{\theta\}) \   \text{and}\ \textbf{A}_0 (M_\theta,\{\theta\})   \    \text{hold with} \\ 
 \hspace{7.01cm} \sum\limits_{k=1}^{\infty} \max \left\{\alpha_k(g),\, \alpha^{(0)}_{k,Y}(f_\theta,\{\theta\}) + \|\xi_0\|_r \alpha^{(0)}_{k,Y}(M_\theta,\{\theta\}\right\} 
 <1
\Big\}\\
\bigcup 
\Big\{
 \theta  \in \R^d \ \big / \ f_\theta=0 \text{ and } \textbf{A}_0 (H_\theta,\{\theta\})  \text{ holds with } 
 \|\xi_0\|^2_r  \sum\limits_{k=1}^{\infty} \max \left\{\alpha_k(g),\,\alpha^{(0)}_{k,Y}(H_\theta,\{\theta\}) \right\}
 <1
\Big\},
\end{multline*}
with the convention (set throughout this section) that if \textbf{A}$_i(M_\theta,\Theta)$ holds, then $\alpha^{(i)}_{k,Y}(H_\theta,\Theta)=\alpha^{(i)}_{k,X}(H_\theta,\Theta)=0$ for all $k\in\N$ and if 
\textbf{A}$_i(H_\theta,\Theta)$ holds then $\alpha^{(i)}_{k,Y}(M_\theta,\Theta)=\alpha^{(i)}_{k,X}(M_\theta,\Theta)=0$ for all $k\in\N$.
\medskip

\noindent 
The following regularity assumptions are also considered in Diop and Kengne (2022a) to establish the consistency and to derive the asymptotic distribution of $\widehat{\theta}(m^*)$.

\medskip
\noindent
  ($\mathcal{AC}.\textbf{A0}$): For all  $\theta \in \Theta$ and some $t \in \Z$, 
 $ \big( f^t_{\theta^*}= f^t_{\theta}  \ \text{and} \ H^t_{\theta^*}= H^t_{\theta} \ \ a.s. \big) \Rightarrow ~ \theta= \theta^*$.
 
\medskip
\noindent
  ($\mathcal{AC}.\textbf{A1}$): $\exists  \underline{h}>0$ such that $\displaystyle \inf_{ \theta \in \Theta} H_\theta(y,x)  \geq \underline{h}$,  for all $(y,x) \in \R^{\infty}\times (\R^{d_x})^{\infty}$.
  
  \medskip
  \noindent
  ($\mathcal{AC}.\textbf{A2}$): 
  \begin{equation} \label{AC_X_Riemann}
   \alpha^{(i)}_{k,Y}(f_\theta,\Theta)+\alpha^{(i)}_{k,X}(f_\theta,\Theta)+\alpha^{(i)}_{k_\theta,Y}(M_\theta,\Theta)+\alpha^{(i)}_{k,X}(M_\theta,\Theta)
   +\alpha^{(i)}_{k_\theta,Y}(H_\theta,\Theta)+\alpha^{(i)}_{k,X}(H_\theta,\Theta) = O(k^{-\gamma}), 
\end{equation} 
 for $i=0,1,2$ and some $\gamma >3/2$  or
 \begin{multline}\label{Assump_sum_coef_AC_X_log_log}
  \sum_{k \geq 1} \frac{1}{\sqrt{k \log \log k}}\sum_{j \geq k} \sum_{i = 0}^{1}
  \big\{ \alpha^{(i)}_{j,Y}(f_\theta,\Theta)+\alpha^{(i)}_{j,X}(f_\theta,\Theta)+\alpha^{(i)}_{j,Y}(M_\theta,\Theta)+\alpha^{(i)}_{j,X}(M_\theta,\Theta)\\
   +\alpha^{(i)}_{j,Y}(H_\theta,\Theta)+\alpha^{(i)}_{j,X}(H_\theta,\Theta)
   \big\} <
   \infty.
\end{multline}

  \medskip
  \noindent
  ($\mathcal{AC}.\textbf{A3}$): for all  $\theta \in \Theta$, $c \in \R^d$, $\Big(c^T \frac{\partial}{\partial \theta}f^0_{\theta^*}=0$ or $c^T \frac{\partial}{\partial    \theta}H^0_{\theta^*}=0 \Big) ~a.s.$   $\Longrightarrow~ c=0$.

\medskip

 Let us check the assumptions (\textbf{A1})-(\textbf{A4}) for this class of models.

\medskip
\noindent
 (i) If $\theta^* \in \Theta \cap \Theta(r)$ with $r\geq 1$, then, there exists a $\tau$-weakly dependent stationary, ergodic and non anticipative solution $(Z_t)_{t\in \Z}$ $Z_t=(Y_t,X_t)$, of (\ref{Model_AC_X}) belonging to $\mathcal {AC}$-$X(M_{\theta^*},f_{\theta^*})$ and satisfying $ \E[{\| Z_0\|}^r] <\infty$ (see Diop and Kengne (2022a)). Which shows that (\textbf{A1}) holds.

\medskip

\noindent (ii) Let us define for all $\theta \in \Theta$   
 \begin{equation}\label{def_fi.t_AC_AC_X}
\varphi_t(\theta):= \dfrac{(Y_t-  f^t_\theta)^2}{ H^t_\theta} +\log(H^t_\theta)
 \end{equation}
  with $f^t_\theta := f_\theta(Y_{t-1},\ldots; X_{t-1},\ldots)$,   
  $M^t_\theta := M_\theta(Y_{t-1},\ldots ; X_{t-1},\ldots)$ and $H^t_\theta := (M^t_\theta)^2$.
  Assume \textbf{A}$_0(f_\theta,\Theta)$, \textbf{A}$_0(M_\theta,\Theta)$ (or \textbf{A}$_0(H_\theta,\Theta)$), ($\mathcal{AC}.\textbf{A0}$), ($\mathcal{AC}.\textbf{A1}$), ($\mathcal{AC}.\textbf{A2}$) ( (\ref{AC_X_Riemann}) with $i=0$) hold and $\theta^* \in \Theta \cap \Theta(2)$.   
  For any $t \in \Z$, the function $\theta \rightarrow \varphi_t(\theta)$  is continuous on $\Theta$ and since $(Y_t,X_t)_{t \in \Z}$ is stationary, ergodic and  $f_\theta (\cdot)$, $M_\theta (\cdot)$ are measurable functions for any $\theta \in \Theta$, then $( \varphi_t(\theta))_{t \in \Z}$ is stationary and ergodic.
 Also, from Lemma 1 and the proof of Theorem 2.1 in Diop and Kengne (2022a), it holds that
 \[ \E \|\varphi_t(\theta) \|_{\Theta} < \infty ,  \text{ and }
     \dfrac{1}{n}  \sum_{t=1}^n \| \widehat{\varphi}_t(\theta) 
     - \varphi_t(\theta)  \|_{\Theta} 
     = o(1) ~ a.s.;\]
 which shows that (\ref{eq_assump_app}) is satisfied.
   The uniqueness of $\theta^*$ satisfying $\theta^* = \underset{\theta \in \Theta}{\text{argmin}}( \E[\varphi_0(\theta)]) $ is established in the proof of Theorem 2.1 in Diop and Kengne (2022a). In the same way, one can show for any $m \in \mathcal{M}$, the existence of a unique $\theta^*(m) \in \Theta(m)$ satisfying $\theta^*(m) = \underset{\theta \in \Theta(m)}{\text{argmin}}( \E[\varphi_0(\theta)]) $. Thus, (\textbf{A2}) is satisfied.

\medskip

\noindent (iii)  Assume \textbf{A}$_i(f_\theta,\Theta)$, \textbf{A}$_i(M_\theta,\Theta)$ (or \textbf{A}$_i(H_\theta,\Theta)$) for $i=0,1$, ($\mathcal{AC}.\textbf{A0}$), ($\mathcal{AC}.\textbf{A1}$), ($\mathcal{AC}.\textbf{A2}$) ((\ref{AC_X_Riemann}) with $i=0,1$) hold and $\theta^* \in \Theta \cap \Theta(4)$. 
Similarly as in (ii), the process $(\partial \varphi_t(\theta)/ \partial \theta)_{t \in \Z}$ is  stationary and ergodic for any $\theta \in \Theta$.
By going as in the proof of Theorem 2.2 in Diop and Kengne (2022a), we get $\E \| \partial \varphi_t(\theta)/ \partial \theta \|^2_{\Theta(m)} < \infty$ and the matrix $G(m)$ exists for any $m \in \mathcal{M}$.
If ($\mathcal{AC}.\textbf{A2}$) ((\ref{AC_X_Riemann}) with $i=0,1$) or ($\mathcal{AC}.\textbf{A2}$) (with (\ref{Assump_sum_coef_AC_X_log_log})) holds then, (\ref{eq_assump_app_dif}) or (\ref{eq_assump_app_dif_log2}) is satisfied respectively, and $\big( \left(  \partial  \varphi_t(\theta^*) / \partial \theta_i \right)_{i \in m},\mathcal{F}_{t}\big)_{t \in \mathbb{Z}}$ is a martingale difference sequence, where $\mathcal{F}_{t}=\sigma((Y_s,X_s),\, s \leq t)$ is the $\sigma$-field generated by the whole past at time $t$; see Lemma 2, 3 and 5 in Diop and Kengne (2022a). 
One can also see that, if $ m^* \subseteq m$, then $G(m)$ is positive definite under ($\mathcal{AC}.\textbf{A3}$). Therefore, (\textbf{A3}) holds. 

\medskip

\noindent (iv) Similarly, from the Lemma 2 in Diop and Kengne (2022a) and the same arguments as in the proof of Theorem 2 in Bardet and Wintenberger (2009), we get that (\textbf{A4}) is satisfied under \textbf{A}$_i(f_\theta,\Theta)$, \textbf{A}$_i(M_\theta,\Theta)$ (or \textbf{A}$_i(H_\theta,\Theta)$) for $i=0,1,2$, ($\mathcal{AC}.\textbf{A0}$), ($\mathcal{AC}.\textbf{A1}$), ($\mathcal{AC}.\textbf{A2}$) ((\ref{AC_X_Riemann}) with $i=0,1,2$) and $\theta^* \in \Theta \cap \Theta(4)$. 

\medskip

\noindent Thus, the results of Theorem \ref{th1_mce} and \ref{th2_mce} apply to the class $\mathcal {AC}$-$X(M_\theta,f_\theta)$.

 \begin{rmrk} \label{rem_AC_X}
   The question considered in this section has been addressed by Diop and Kengne (2022a); see also Bardet {\it et al.} (2020), Kengne (2021), Kamila (2021) for affine causal models without covariates.
  \begin{enumerate}
  \item In order to establish the weak consistency, these authors imposed $\kappa_n/\sqrt{\log \log n} \limiten \infty$ or a dependent relation between the regularization parameter $\kappa_n$ and the Lipschitz coefficients of $f_\theta$ and $M_\theta$ (or $H_\theta$). The results obtained here show that such conditions is not needed; we only need that $\kappa_n$ tends to infinity, the rate does not matter. 
  \item The results obtained here go beyond those of these authors, by providing a class of $(\kappa_n)$ that does not ensure the consistency. 
   \item These results for the class $\mathcal {AC}$-$X(M_\theta,f_\theta)$ generalize those of Hannan and Deistler (2012) (Theorem 5.4.1) established for ARMAX model.
  \end{enumerate}
 \end{rmrk}

\section{Application to a class of observation-driven multivariate count time series}\label{Sec_class_OD}
%
 Let $\{Y_{t}= (Y_{t,1},\ldots,Y_{t,d_{y}})^T,\,t\in \Z \}$ a multivariate count time series  with value in $\N_0^{d_{y}}$ (with $d_{y} \in \N$, $\N_0 = \N \cup \{ 0 \}$) and denote by $\mathcal{F}_{t}=\sigma\left\{Y_{t},\ldots  \right\}$ the $\sigma$-field generated by the whole past at time $t$.
 Consider the class of observation-driven multivariate count time series defined by (see Diop and Kengne (2021))
 
 \medskip
 \noindent 
 {\bf Class} $\mathcal{MOD}(f_{\theta})$: A multivariate count process $Y=\{Y_{t},\,t\in \mathcal \Z  \}$ belongs to $\mathcal{MOD}(f_{\theta})$ if it satisfies:
  \begin{equation} \label{Model_MOD}
             \E(Y_t|\mathcal{F}_{t-1})=f_{\theta}(Y_{t-1},Y_{t-2},\ldots) ~~ \forall t \in \Z ,
   \end{equation}
where $f_{\theta}(\cdot)$ is a measurable multivariate function with  non-negative components, assumed to be known up to the  parameter $\theta$, which belongs to a compact subset $ \Theta\subset \R^d$ ($d \in \N$).
   This class includes numerous classical univariate integer-valued time series, such as; the Poisson INGARCH models (proposed by  Ferland {\it et al.} (2006)), the negative binomial INGARCH models (see Zhu (2011)), the binomial INGARCH (see for instance Wei$ß$ and Pollett (2014)), the Poisson exponential autoregressive models (proposed by Fokianos {\it et al.} (2009)), among others. 
  Multivariate count time series such as, the bivariate Poisson INGARCH model (see Lee  \textit{et al.} (2018)), the flexible bivariate Poisson INGARCH model (see Cui {\it et al.} (2019)), the multivariate count autoregression (see Fokianos {\it et al.} (2020)) are specific examples that belong to $\mathcal{MOD}(f_{\theta})$. 
 
 \medskip
 
  Consider a trajectory $(Y_1,\ldots,Y_n)$ generated from a process $(Y_t)_{t \in \Z}$, assumed to be stationary and ergodic (that is, (\textbf{A1}) is satisfied) and belongs to $\mathcal{MOD}(f_{\theta^*})$, where $\theta^* \in \Theta \subset \R^d$ is the true parameter.    
  Let us carry out the model selection problem presented in Section \ref{sec_intro}, with the finite collection of the competing model  $\mathcal M$  that contains the true model $m^*$ corresponding to the parameter $\theta^*$. 
In this section, assume that 
 \begin{equation}\label{moment_MOD}
    \exists C>0, \epsilon >0, \text{ such that } \forall t \in \Z, ~ ~ \|Y_{t}\|_{1+\epsilon} <C. 
   \end{equation}
   
   \medskip  
   
 Define the contrast function based on the conditional Poisson quasi log-likelihood given (up to a constant) by
\begin{equation*}
\widehat{\Phi}_n( \theta) \coloneqq  \sum_{t=1}^{n} \widehat{\varphi}_t(\theta) ~ \text{ with }~
 \widehat{\varphi}_t(\theta) =  - \sum_{i=1}^{d_{y}} \left(Y_{t,i} \log \widehat{\lambda}_{t,i}(\theta)- \widehat{\lambda}_{t,i}(\theta)\right),
 \end{equation*}
  where $ \widehat{\lambda}_t(\theta) \coloneqq \big(\widehat{\lambda}_{t,1}(\theta), \ldots, \widehat{\lambda}_{t,d_{y}}(\theta) \big)^T =  f_\theta(Y_{t-1}, \ldots, Y_{1},0,\ldots)$.
Thus, the MCE given a model $m$ is defined by
 \begin{equation*}
  \widehat{\theta} (m) \coloneqq  \underset{\theta\in \Theta(m)}{\text{argmin}} \big(\widehat{\Phi}_n(\theta)\big).
  \end{equation*}

  \medskip
  
  \noindent
  For a process $\{Y_t,\, t \in \Z\}$ of the class $\mathcal{MOD}(f_{\theta^{*}})$, Diop and Kengne (2021) have set the following assumptions in order to study the stability properties when the distribution of $Y_t|\mathcal{F}_{t-1}$ belongs to the $d_y$-parameter exponential family  and to establish the consistency and the asymptotic normality of $\widehat{\theta} (m^*)$. 

   \medskip
   \noindent
  \textbf{Assumption} \textbf{A}$_i (\Theta)$ ($i=0,1,2$):
    For any $y \in \big(\mathbb{N}_0^{d_y} \big)^{\infty}$, the function $\theta \mapsto f_\theta(y)$ is $i$ times continuously differentiable on $\Theta$  with $ \left\| \partial^i f_\theta(0)/ \partial \theta^i\right\|_\Theta<\infty $; 
    and
      there exists a sequence of non-negative real numbers $(\alpha^{(i)}_k)_{k\geq 1} $ satisfying
     $ \sum\limits_{k=1}^{\infty} \alpha^{(0)}_k <1 $ (or $ \sum\limits_{k=1}^{\infty} \alpha^{(i)}_k <\infty $ for $i=1, 2$);
   such that for any  $y, y' \in \big(\mathbb{N}_0^{m} \big)^{\infty}$,
  \[ 
   \Big \| \frac{\partial^i f_\theta(y)}{ \partial \theta^i}-\frac{\partial^i f_\theta(y')}{\partial\theta^i} \Big \|_\Theta
  \leq  \sum\limits_{k=1}^{\infty}\alpha^{(i)}_k \|y_k-y'_k\|.
  \]

  \medskip
  \noindent
 ($\mathcal{MOD}.\textbf{A0}$): 
     For all  $\theta \in \Theta$,
 $ \big( f_{\theta^*}(Y_{t-1}, Y_{t-2}, \ldots) \equalpsn  f_{\theta}(Y_{t-1}, Y_{t-2}, \ldots)  ~ \text{ for some } t \in \Z \big) \Rightarrow ~ \theta^* = \theta$; 
 moreover, $\exists  \underline{c}>0$ such that $ f_\theta(y)  \geq \underline{c} \textbf{1}_{d_y}$ componentwise, for all $\theta \in \Theta$, $ y \in  \big( \N_0^{d_y} \big)^{\infty} $, where $\textbf{1}^T_{d_y}=(1,\ldots,1)$ is a vector of dimension $d_y$.
%

 %
 
 \medskip
 \noindent
  ($\mathcal{MOD}.\textbf{A1}$): 
\begin{equation}\label{Assum_Riemann_MOD}
 \alpha^{(i)}_{k}= \mathcal{O}(k^{-\gamma}) 
\end{equation}
for $i=0,1,2$ and  some $\gamma >3/2$ or
 \begin{equation}\label{Assump_sum_coef_MOD_log_log}
  \sum_{k \geq 1} \frac{1}{\sqrt{k \log \log k}}\sum_{j \geq k} 
  \big( \alpha^{(0)}_{j} + \alpha^{(1)}_{j} \big) < \infty.
\end{equation}
 
 \medskip
 \noindent
  ($\mathcal{MOD}.\textbf{A2}$):  The family $ \big(\frac{\partial \lambda_{t} (\theta^* )}{\partial    \theta_i} \big)_{1\leq i \leq d} $  is $a.e.$  linearly independent.

 \medskip
 
  \medskip
   
 According to Diop and Kengne (2021), it holds that, under \textbf{A}$_0(\Theta)$, there exists a $\tau-weakly$ dependent, stationary and ergodic process $\{ Y_t, ~ t \in \Z \}$ satisfying $\E\|Y_t \| < \infty$ and solution of
 \begin{equation*} 
  Y_t|\mathcal{F}_{t-1} \sim p(y|\eta_t) ~ \text{ with } ~ \lambda_t(\theta) := \E(Y_t|\mathcal{F}_{t-1})=f_{\theta}(Y_{t-1},Y_{t-2},\ldots)  
 \end{equation*}
 where $p(\cdot | \cdot )$ is a multivariate discrete distribution belonging to the $d_y$-parameter exponential family and  $\eta_t$ is the natural parameter of the distribution of $Y_t|\mathcal{F}_{t-1}$. Thus, (\textbf{A1}) holds for such distribution family. Also, these authors have established the consistency and the asymptotic normality of $\widehat{\theta} (m^*)$. From their Lemma 6.2 and 6.3 and the proofs, and with the same arguments as in Section \ref{Sec_class_AC} with 
 $ \varphi_t(\theta) =  - \sum_{i=1}^{d_{y}} \left(Y_{t,i} \log \lambda_{t,i}(\theta)- \lambda_{t,i}(\theta)\right)$,
  where $ \lambda_t(\theta) \coloneqq \big(\lambda_{t,1}(\theta), \ldots, \lambda_{t,d_{y}}(\theta) \big)^T =  f_\theta(Y_{t-1}, Y_{t-2},\ldots)$;  one gets that (\textbf{A2})-(\textbf{A4}) are satisfied under \textbf{A}$_i (\Theta)$ ($i=0,1,2$) and ($\mathcal{MOD}.\textbf{A0}$)-($\mathcal{MOD}.\textbf{A2}$).   
 
 \medskip

\noindent Therefore, the results of Theorem \ref{th1_mce} and \ref{th2_mce} apply to the class $\mathcal{MOD}(f_{\theta})$.

 \begin{rmrk} \label{rem_MOD}
  Diop and Kengne (2022b) have considered the case of univariate observation-driven integer-valued time series and establish the weak consistency of the model selection procedure based on the Poisson quasi likelihood. Besides dealing with multivariate models here, we establish the strong consistency, address the non consistency issue and the conditions on the regularization parameter $\kappa_n$ is weaker.  
 \end{rmrk}








 \section{Proofs of the main results}  
  Throughout the sequel, $C$ denotes a positive constant whom value may differ from an inequality to another; 
  $O(\cdot)$, $o(\cdot)$ is always meant that the indicated order relation holds $a.s.$.
 
\subsection{Proof of Theorem \ref{th1_mce}}
 \noindent 1. It suffices to show that 
\begin{equation}\label{cond_proof_th1_mce}
\lim\limits_{n \rightarrow \infty} P (\widehat m \supsetneq  m^*) = \lim\limits_{n \rightarrow \infty}P (\widehat m \nsupseteq m^*)=0.
\end{equation} 
%
\indent (i) Let a model $m\in {\cal M}$ such as $m  \supsetneq  m^*$.
 In the proof of this first part, we set $\partial \psi (\theta) / \partial \theta := (\partial \psi(\theta) / \partial \theta_i)_{i \in m}$ for any function $\psi$ and $\theta \in \Theta$.
We have, 
\begin{equation}\label{C_m_star_C_m}
  \widehat{C}(m) - \widehat{C}(m^*)  =  \widehat{\Phi}_n\big(\widehat{\theta}(m)\big)-\widehat{\Phi}_n\big(\widehat{\theta}(m^*) \big) + \kappa_n  ( |m|-|m^*|). 
\end{equation}
Let us prove that
\begin{equation}\label{Phi_m_star_Phi_m}
 \widehat{\Phi}_n\big(\widehat{\theta}(m)\big)-\widehat{\Phi}_n\big(\widehat{\theta}(m^*)\big)  = O_P(1) .
\end{equation}  
 We have $\theta^* \in \Theta(m) \cap \interior{\Theta}$ and $ \widehat{\theta}(m) \limitepsn \theta^*$ (see Remark \ref{rm_assump}), therefore, for $n$ large enough 
 $ \dfrac{\partial\widehat{\Phi}_n(\widehat{\theta}(m) )}{\partial \theta} =0 $.  
Hence, from the Taylor expansion of $\widehat{\Phi}_n$, we can find $\overline{\theta}(m)$ between $\widehat{\theta}(m)$ and $\theta^*$ such that
\begin{equation}\label{Taylor_widehat_Phi}
\widehat{\Phi}_n\big(\theta^*) - \widehat{\Phi}_n\big(\widehat{\theta}(m)\big) =  \frac{1}{2}\big( \widehat{\theta}(m) - \theta^* \big)^T \dfrac{\partial^2 \widehat{\Phi}_n\big( \overline{\theta}(m) \big)}{\partial \theta^2} \big( \widehat{\theta}(m) - \theta^* \big).
\end{equation}
Moreover, for any $i \in m$, we can find $\dot{\theta}_i(m)$ between $\widehat{\theta}(m)$ and $\theta^*$ such that, for $n$ large enough,
 \[ 0 =  \dfrac{\partial\widehat{\Phi}_n(\widehat{\theta}(m) )}{\partial \theta_i} = \dfrac{\partial\widehat{\Phi}_n(\theta^*)}{\partial \theta_i} +   \dfrac{\partial^2 \widehat{\Phi}_n\big( \dot{\theta}_i(m) \big)}{\partial \theta \partial \theta_i} \big( \widehat{\theta}(m) - \theta^* \big) .\]
Therefore, we obtain for $n$ large enough, 
\begin{equation}\label{theta_m_theta_star_partial_widehat_Phi}
\widehat{\theta}(m) - \theta^* = - \dfrac{1}{n} \widehat{F}_n^{-1}(m) \dfrac{\partial\widehat{\Phi}_n(\theta^*)}{\partial \theta} ~ ~ \text{ where } ~ 
 \widehat{F}_n(m) = \dfrac{1}{n} \Big( \dfrac{\partial^2 \widehat{\Phi}_n\big( \dot{\theta}_i(m) \big)}{\partial \theta \partial \theta_i} \Big)_{i \in m}. 
\end{equation}
For this case of overfitting, we have for $i \in m$, $\widehat{\theta}(m),~ \overline{\theta}(m),~ \dot{\theta}_i(m) \limitepsn \theta^*$.
In addition to (\textbf{A4}), it holds that
\begin{equation}\label{F_n_partial2_widehat_Phi_conv}
\widehat{F}_n(m) \limitepsn F(m)  \text{ and }  \dfrac{1}{n} \dfrac{\partial^2 \widehat{\Phi}_n\big( \overline{\theta}(m) \big)}{\partial \theta^2} \limitepsn F(m) ~  \text{ where } ~ F(m)= \Big(\E \Big[ \dfrac{ \partial^2 \varphi_0(\theta^*)}{ \partial\theta_i \partial\theta_j}  \Big] \Big)_{i,j \in m}  
\end{equation}     
is invertible in this case.
Thus, for $n$ large enough and with a sufficiently large probability, $\widehat{F}_n(m)$ is invertible.
Hence, according to (\ref{Taylor_widehat_Phi}), (\ref{theta_m_theta_star_partial_widehat_Phi}), (\ref{F_n_partial2_widehat_Phi_conv}) and (\textbf{A3}) with (\ref{eq_assump_app_dif}), it holds that
\begin{align}
\nonumber \widehat{\Phi}_n\big(\widehat{\theta}(m)\big)-\widehat{\Phi}_n\big(\theta^*) &= - \frac{1}{2 n^2} \dfrac{\partial\widehat{\Phi}_n(\theta^*)}{\partial \theta^T}  \widehat{F}_n^{-1}(m)  \dfrac{\partial^2 \widehat{\Phi}_n\big( \overline{\theta}(m) \big)}{\partial \theta^2} \widehat{F}_n^{-1}(m) \dfrac{\partial\widehat{\Phi}_n(\theta^*)}{\partial \theta} \\
\label{Phi_m_theta_star_partial_F_n} 
 &= - \frac{1}{2} \Big( \frac{1}{\sqrt{n}} \dfrac{\partial\widehat{\Phi}_n(\theta^*)}{\partial \theta^T} \Big)  \widehat{F}_n^{-1}(m) \Big( \dfrac{1}{n} \dfrac{\partial^2 \widehat{\Phi}_n\big( \overline{\theta}(m) \big)}{\partial \theta^2} \Big) \widehat{F}_n^{-1}(m) \Big( \frac{1}{\sqrt{n}} \dfrac{\partial\widehat{\Phi}_n(\theta^*)}{\partial \theta} \Big) \\
\label{Phi_m_theta_star_partial}
&=  - \Big( \frac{1}{\sqrt{n}} \dfrac{\partial \Phi_n(\theta^*)}{\partial \theta^T}  + o_P(1) \Big) O(1) \Big( \frac{1}{\sqrt{n}} \dfrac{\partial \Phi_n(\theta^*)}{\partial \theta}  + o_P(1) \Big).
\end{align}
According to (\textbf{A3}) and the central limit theorem for stationary ergodic martingale difference sequence, it holds that
\begin{equation} \label{prov_th1_mce_CLTMDS}
 \frac{1}{\sqrt{n}} \dfrac{\partial \Phi_n(\theta^*)}{\partial \theta} = O_P(1).
\end{equation}
Therefore, (\ref{Phi_m_theta_star_partial}) gives
\begin{equation}\label{Phi_m_theta_star_O_P_1}
 \widehat{\Phi}_n\big(\widehat{\theta}(m)\big)-\widehat{\Phi}_n\big(\theta^*) = O_P(1). 
\end{equation}
By using the same arguments with $m=m^*$, it follows that 
\begin{equation}\label{Phi_m_star_theta_star_O_P_1}
 \widehat{\Phi}_n\big(\widehat{\theta}(m^*)\big)-\widehat{\Phi}_n\big(\theta^*) = O_P(1). 
\end{equation}
Hence, (\ref{Phi_m_star_Phi_m}) follows from (\ref{Phi_m_theta_star_O_P_1}) and (\ref{Phi_m_star_theta_star_O_P_1}).
Thus, according to (\ref{C_m_star_C_m}), (\ref{Phi_m_star_Phi_m}) and  since $|m^*| < |m|$, we get
\begin{equation*} 
  \widehat{C}(m) - \widehat{C}(m^*) \limiteproban \infty.  
\end{equation*}  
This implies,
\begin{equation*} 
  \lim_{n \rightarrow \infty} P\big( \widehat{C}(m) - \widehat{C}(m^*) > 0\big) =1,  
\end{equation*} 
when $m \supsetneq m^*$.
Since $ \widehat{C}( \widehat m ) - \widehat{C}(m^*) \leq 0 ~ a.s.$, we get
\begin{equation}\label{m_star_m_hat_proba}
  P (\widehat m \supsetneq  m^*) \limiten 0 .  
\end{equation}  

\medskip

\indent (ii) Let $m\in {\cal M}$ such as $m  \nsupseteq  m^*$. 
  We have,
   \begin{equation}\label{C_m_star_C_m_n}
  \widehat{C}(m) - \widehat{C}(m^*)  =  \widehat{\Phi}_n\big(\widehat{\theta}(m)\big)-\widehat{\Phi}_n\big(\widehat{\theta}(m^*)\big)  +  \kappa_n ( |m|-|m^*|). 
\end{equation}
Also,
\begin{multline}\label{Phi_n_theta_m_star_m}
\widehat \Phi_n \big(\widehat  \theta(m) \big)-\widehat \Phi_n \big(\widehat  \theta(m^*) \big) 
= \left(\widehat \Phi_n \big(\widehat  \theta(m) \big)  - \Phi_n \big(\widehat  \theta(m) \big) \right) -\left(\widehat \Phi_n \big(\widehat  \theta(m^*) \big)  - \Phi_n \big(\widehat  \theta(m^*) \big) \right)\\
+ \left( \Phi_n \big(\widehat  \theta(m) \big)  - \Phi_n \big(\widehat  \theta(m^*) \big) \right).
\end{multline}
From the assumption (\textbf{A2}), we have, 
\begin{equation}\label{Phi_n_Phi_o_n}
\widehat \Phi_n \big(\widehat  \theta(m) \big)  - \Phi_n \big(\widehat  \theta(m) \big)  =  \widehat \Phi_n \big(\widehat  \theta(m^*) \big)  - \Phi_n \big(\widehat  \theta(m^*) \big)  = o(n),
\end{equation} 
and in addition to the uniform strong law of large number applied to $\big(\varphi_t(\theta) \big)_{t \in \Z}$ for $\theta \in \Theta$, we get
\begin{equation} \label{Phi_n_Phi_o_n_def_phi}
\dfrac{1}{n} \Phi_n \big(\widehat  \theta(m) \big)  - \phi \big(\widehat  \theta(m) \big)  =  \dfrac{1}{n} \Phi_n \big(\widehat  \theta(m^*) \big)  - \phi \big(\widehat  \theta(m^*) \big)  = o(1), \text{ where } 
\phi(\theta) = \E [ \varphi_0(\theta) ] \text{ for all } \theta \in \Theta. 
\end{equation}
Hence,
\begin{align}\label{Phi_n_theta_m_star_m_phi}
\nonumber \dfrac{1}{n} \Big( \Phi_n \big(\widehat  \theta(m) \big)  - \Phi_n \big(\widehat  \theta(m^*) \big) \Big) &=  \phi \big(\widehat  \theta(m) \big)  - \phi \big(\widehat  \theta(m^*) \big)  + o(1) \\
&= \big(\phi \big(\widehat  \theta(m) \big) - \phi \big( \theta^*(m) \big) - \big( \phi \big(\widehat  \theta(m^*) \big) - \phi \big(\theta^* \big) \big) + \big( \phi(\theta^*(m)) - \phi( \theta^*) \big) +  o(1) .
\end{align}
According to the consistency of  $\widehat  \theta(m)$ and $ \theta^*(m)$ (see Remark \ref{rm_assump}), it holds that 
$  \phi \big(\widehat  \theta(m) \big) - \phi \big( \theta^*(m) = \phi \big(\widehat  \theta(m^*) \big) - \phi \big(\theta^* \big)  = o(1)$.
Therefore, from (\ref{C_m_star_C_m_n}), (\ref{Phi_n_theta_m_star_m}), (\ref{Phi_n_Phi_o_n}) and (\ref{Phi_n_theta_m_star_m_phi}), we have
  \begin{align} \label{prov_th1_mce_approx_Cm_Cm_start}
 \frac{1}{  n} \big( \widehat{C}(m) - \widehat{C}(m^*) \big) = \phi(\theta^*(m)) - \phi( \theta^*) + \frac{\kappa_n}{ n}  ( |m|-|m^*|) + o(1). 
\end{align} 
 By virtue of (\textbf{A2}), the function $\phi : \Theta \rightarrow \R$ has a unique minimum at $\theta^*$. In this case where $m \nsupseteq m^*$, we have $\theta^* \notin \Theta(m)$.
 As a consequence, $\phi(\theta^*(m)) - \phi( \theta^*) >0$. Since $\kappa_n = o(n)$, it holds from (\ref{prov_th1_mce_approx_Cm_Cm_start}) that, $\widehat{C}(m) - \widehat{C}(m^*)  > 0 ~ a.s.$ for $n$ large enough.
 This implies (since $\widehat{C}( \widehat m ) - \widehat{C}(m^*) \leq 0~ a.s.$) $P(\widehat{m}  \nsupseteq  m^*) \limiten 0$.
Hence, in addition to (\ref{m_star_m_hat_proba}), (\ref{cond_proof_th1_mce}) holds. This completes the proof of the first part of the theorem.
 
\medskip

\medskip

Before proving the second part of the theorem, let us consider the following lemma.
\begin{lem}\label{lem_a_n_theta_m_theta_star}
Assume that the conditions of Theorem \ref{th1_mce} hold. 
For any $m \in \mathcal{M}$ such that $m^* \subseteq m$, it holds that,
\begin{equation}\label{a_n_theta_m_tilde_theta_star}
\sqrt{n}
\begin{pmatrix}
    \Big( \big(\widehat{\theta}(m) \big)_i - (\theta^*)_i \Big)_{i \in m^*}   \\
    \Big( \big(\widehat{\theta}(m^*) \big)_i - (\theta^*)_i \Big)_{i \in m} 
  \end{pmatrix}
\limiteloin 
 \mathcal{N}\big( 0 , \Sigma(m^*,m) \big), 
\end{equation}
where
\[ 
 \Sigma(m^*,m) = \begin{pmatrix}
    F(m^*)^{-1} G(m^*) F(m^*)^{-1} & F(m^*)^{-1} G(m^*, m) F(m)^{-1}\\
   F(m)^{-1} G( m, m^*) F(m^*)^{-1}  & F(m)^{-1} G(m) F(m)^{-1}
  \end{pmatrix},
  \]
and  
$G(m^*, m) = G(m, m^*)^T  = \E \left[ \Big( \dfrac{\partial \varphi_0(\theta^*)}{\partial \theta_i} \Big)_{i \in m^*} 
\Big( \dfrac{\partial \varphi_0(\theta^*)}{\partial \theta_i} \Big)_{i \in m}^T \right]$.
\end{lem}

\medskip

\medskip

\noindent 2. Assume $ \displaystyle{\limsup_{n \rightarrow \infty} \kappa_n < \infty}$ and that there exists $\widetilde{m} \in \mathcal{M}$ such that $m^* \subsetneq \widetilde{m} $.
 From (\ref{Taylor_widehat_Phi}), (\ref{F_n_partial2_widehat_Phi_conv}) and Lemma \ref{lem_a_n_theta_m_theta_star}, we get
 \begin{equation*} 
 \widehat{\Phi}_n\big(\theta^*) - \widehat{\Phi}_n\big(\widehat{\theta}(\widetilde{m})\big) =  \frac{n}{2}\big( \widehat{\theta}(\widetilde{m}) - \theta^* \big)^T F(\widetilde{m}) \big( \widehat{\theta}(\widetilde{m}) - \theta^* \big) + o_P(1).
\end{equation*}
The same arguments with $m^*$ yield
 \begin{equation*} 
 \widehat{\Phi}_n\big(\theta^*) - \widehat{\Phi}_n\big(\widehat{\theta}(m^*)\big) =  \frac{n}{2}\big( \widehat{\theta}(m^*) - \theta^* \big)^T F(m^*) \big( \widehat{\theta}(m^*) - \theta^* \big) + o_P(1).
\end{equation*}
 Hence, according to Lemma \ref{lem_a_n_theta_m_theta_star} and Lemma 3.2 in Vuong (1989), we have
  \begin{align} \label{Taylor_widehat_Phi_m_tilde_m_star_oP}
\nonumber & 2 \Big( \widehat{\Phi}_n\big(\widehat{\theta}(m^*)\big) - \widehat{\Phi}_n\big(\widehat{\theta}(\widetilde{m})\big) \Big)  \\
\nonumber &= n\big( \widehat{\theta}(\widetilde{m}) - \theta^* \big)^T F(\widetilde{m}) \big( \widehat{\theta}(\widetilde{m}) - \theta^* \big) - n\big( \widehat{\theta}(m^*) - \theta^* \big)^T F(m^*) \big( \widehat{\theta}(m^*) - \theta^* \big) + o_P(1) \\
\nonumber  &= 
 n
\begin{pmatrix}
    \Big( \big(\widehat{\theta}(m) \big)_i - (\theta^*)_i \Big)_{i \in m^*}   \\
    \Big( \big(\widehat{\theta}(m^*) \big)_i - (\theta^*)_i \Big)_{i \in m} 
  \end{pmatrix}
  ^T
  Q(m^*,\widetilde{m})
  \begin{pmatrix}
    \Big( \big(\widehat{\theta}(m) \big)_i - (\theta^*)_i \Big)_{i \in m^*}   \\
    \Big( \big(\widehat{\theta}(m^*) \big)_i - (\theta^*)_i \Big)_{i \in m} 
  \end{pmatrix} 
  + o_P(1) \\
 & \limiteloin W(|m^*| + |\widetilde{m}|, \lambda),
  \end{align}
where 
\[Q(m^*,\widetilde{m}) =  
   \begin{pmatrix}
   -F(m^*)  & 0  \\
     0   & F(\widetilde{m})
  \end{pmatrix}
   ,\]
$W(|m^*| + |\widetilde{m}|, \lambda)$ is a weighted sum of chi-squares with parameters $(|m^*| + |m|, \lambda)$ and $\lambda$ is the vector of eigenvalues of $Q(m^*,\widetilde{m}) \Sigma(m^*,\widetilde{m})$.
These eigenvalues are all real (see Lemma 3.2 in Vuong (1989)). 
 Recall that, $W(|m^*| + |\widetilde{m}|, \lambda) = \sum_{j=1}^{|m^*| + |\widetilde{m}|} \lambda_j Z^2_j$, where $\lambda = (\lambda_1,\cdots,\lambda_{|m^*| + |\widetilde{m}|})^T$ and $(Z_j)$ are i.i.d. standard normal random variables.
On can easily see that
\[ Q(m^*,\widetilde{m}) \Sigma(m^*,\widetilde{m})=
    \begin{pmatrix}
   -G(m^*) F(m^*)^{-1}  & -G(m^*,\widetilde{m}) F(\widetilde{m})^{-1}  \\
     G(\widetilde{m}, m^*) F(m^*)^{-1}   &  G(\widetilde{m}) F(\widetilde{m})^{-1}  
  \end{pmatrix}
    . \]
Since $\Sigma(m^*,\widetilde{m})$ is symmetric and positive definite, denote $\Sigma(m^*,\widetilde{m})^{1/2}$ its unique square root. Then, as $Q(m^*,\widetilde{m})$ is non-singular, $\lambda$ is also the eigenvalues of $\Sigma(m^*,\widetilde{m})^{1/2} Q(m^*,\widetilde{m}) \Sigma(m^*,\widetilde{m})^{1/2}$ (see Theorem 1.3.22 in Horn and Johnson (2012)).
Moreover, from the Sylvester's law of inertia, the matrices  $\Sigma(m^*,\widetilde{m})^{1/2} Q(m^*,\widetilde{m}) \Sigma(m^*,\widetilde{m})^{1/2}$ and  $ Q(m^*,\widetilde{m}) $ have the same number of positive, negative and zero eigenvalues.
Thus, there are $|m^*|$ strictly negative and $|\widetilde{m}|$ strictly positive  eigenvalues of $\Sigma(m^*,\widetilde{m})^{1/2} Q(m^*,\widetilde{m}) \Sigma(m^*,\widetilde{m})^{1/2}$ (with $|m^*| < |\widetilde{m}|$). Consequently, the support of a distribution of $W(|m^*| + |\widetilde{m}|, \lambda)$ contains the set $(0,\infty)$. 
Thus, we obtain from (\ref{Taylor_widehat_Phi_m_tilde_m_star_oP}),
\begin{align*}
 \limsup_{n \rightarrow \infty} P\big( \widehat{C}(\widetilde{m}) - \widehat{C}(m^*) < 0  \big) &=  \limsup_{n \rightarrow \infty} P\Big(\widehat{\Phi}_n\big(\widehat{\theta}(m^*) - \widehat{\Phi}_n\big(\widehat{\theta}(\widetilde{m})\big)\big) > \kappa_n  ( |\widetilde{m}|-|m^*|) \Big) \\
 &= P\Big( W(|m^*| + |\widetilde{m}|, \lambda) > 2 \kappa ( |\widetilde{m}|-|m^*|) \Big) >0,
\end{align*}   
where $\kappa = \limsup_{n \rightarrow \infty} \kappa_n$. 
Hence,
\[ \limsup_{n \rightarrow \infty} P( \widehat{m} = m^* ) \leq  \limsup_{n \rightarrow \infty} P\big( \widehat{C}(m^*) \leq  \widehat{C}(\widetilde{m})  \big)  
= \limsup_{n \rightarrow \infty} P\big( \widehat{C}(\widetilde{m}) - \widehat{C}(m^*) \geq 0 \big) < 1,
    \]  
which shows that, the weak consistency fails.

\medskip

\noindent Note that, if  $m^* \not\subset m$ for all $m \in \mathcal{M} \setminus \{m^*\}$, it holds from the proof of of the first part (ii) that $P(\widehat{m}  \nsupseteq  m^*) \limiten 0$. This shows that the consistency holds (according to (\ref{cond_proof_th1_mce})) in this case.

\begin{flushright}
$\blacksquare$ 
\end{flushright}

\noindent 
{\bf Proof of Lemma \ref{lem_a_n_theta_m_theta_star}.}\\
From (\ref{theta_m_theta_star_partial_widehat_Phi}), (\ref{F_n_partial2_widehat_Phi_conv}), (\ref{prov_th1_mce_CLTMDS}) and (\textbf{A3}) with (\ref{eq_assump_app_dif}), we have
\begin{equation*}  
\sqrt{n} \big( \widehat{\theta}(m) - \theta^* \big) = - \big( F(m)^{-1}+ o(1) \big) \left( \dfrac{1}{\sqrt{n}} \Big(\dfrac{\partial \Phi_n(\theta^*)}{\partial \theta_i} \Big)_{i \in m} + o_P(1) \right)
= - \dfrac{1}{\sqrt{n}} F(m)^{-1} \Big( \dfrac{\partial \Phi_n(\theta^*)}{\partial \theta_i} \Big)_{i \in m} + o_P(1).    
\end{equation*}
Similarly, with $m^*$, it follows that
\begin{equation*}  
\sqrt{n} \big( \widehat{\theta}(m^*) - \theta^* \big) = - \big( F(m^*)^{-1}+ o(1) \big) \left( \dfrac{1}{\sqrt{n}} \Big(\dfrac{\partial \Phi_n(\theta^*)}{\partial \theta} \Big)_{i \in m^*} + o_P(1) \right)
= - \dfrac{1}{\sqrt{n}} F(m^*)^{-1} \Big( \dfrac{\partial \Phi_n(\theta^*)}{\partial \theta_i} \Big)_{i \in m^*} + o_P(1)    
\end{equation*}
Hence,
\begin{equation} \label{theta_m_theta_star_F_inv_partial_Phi}
\sqrt{n}
\begin{pmatrix}
    \Big( \big(\widehat{\theta}(m) \big)_i - (\theta^*)_i \Big)_{i \in m^*}   \\
    \Big( \big(\widehat{\theta}(m^*) \big)_i - (\theta^*)_i \Big)_{i \in m} 
  \end{pmatrix}
  = 
 -\dfrac{1}{\sqrt{n}}
\begin{pmatrix}
 F(m)^{-1} \Big(\partial \Phi_n(\theta^*)/\partial \theta_i \Big)_{i \in m^*}  \\
F(m^*)^{-1} \Big(\partial \Phi_n(\theta^*)/\partial \theta_i \Big)_{i \in m} 
\end{pmatrix} 
+ o_P(1) .
\end{equation}
By virtue of central limit theorem for the martingale difference sequence, it holds that,
\begin{equation} \label{CLTMDS_m_m_star} 
\dfrac{1}{\sqrt{n}}
\begin{pmatrix}
 \Big(\partial \Phi_n(\theta^*)/\partial_i \Big)_{i \in m^*}  \\
\Big(\partial \Phi_n(\theta^*)/\partial_i \Big)_{i \in m} 
\end{pmatrix}
\limiteloin \mathcal{N}\left( 0 ,
   \begin{pmatrix}
    G(m^*) & G(m^*, m) \\
    G( m, m^*) & G(m)
  \end{pmatrix}  
  \right).
\end{equation}
Thus, the lemma follows from (\ref{theta_m_theta_star_F_inv_partial_Phi}) and (\ref{CLTMDS_m_m_star}).

\begin{flushright}
$\blacksquare$ 
\end{flushright}

\subsection{Proof of Theorem \ref{th2_mce}}

\noindent 1. 
  Let $m\in {\cal M}$ satisfying $m  \supsetneq  m^*$. 
   In the proof of this first part, we set $\partial \psi (\theta) / \partial \theta := (\partial \psi(\theta) / \partial \theta_i)_{i \in m}$ for any function $\psi$ and $\theta \in \Theta$.
  Recall that 
  \begin{equation} \label{C_m_star_C_m_ps}
   \widehat{C}(m) - \widehat{C}(m^*)  = \widehat{\Phi}_n\big(\widehat{\theta}(m)\big)-\widehat{\phi}_n\big(\widehat{\theta}(m^*)\big)  +  \kappa_n ( |m|-|m^*|). 
  \end{equation}
\noindent Let us establish that, there exists a constant $C>0$ such that,
\begin{equation}\label{Phi_m_star_Phi_m_limsup}
 \limsup_{n \rightarrow \infty}  \dfrac{1}{\log \log n} \big| \widehat{\Phi}_{n} \big(\widehat{\theta}(m)\big)-\widehat{\Phi}_{n}\big(\widehat{\theta}(m^*)\big) \big | \leq C ~ a.s.
\end{equation}  
 From (\ref{Phi_m_theta_star_partial_F_n}) and (\textbf{A3}) with (\ref{eq_assump_app_dif_log2}), we get
 \begin{align}
 \nonumber \big| \widehat{\Phi}_n\big(\widehat{\theta}(m)\big)-\widehat{\Phi}_n\big(\theta^*) \big| & \leq \frac{1}{2} \Big \| \frac{1}{\sqrt{n}} \dfrac{\partial\widehat{\Phi}_n(\theta^*)}{\partial \theta} \Big \|^2  \| \widehat{F}_n^{-1}(m) \|^2 \Big\| \dfrac{1}{n} \dfrac{\partial^2 \widehat{\Phi}_n\big( \overline{\theta}(m) \big)}{\partial \theta^2} \Big \|    \\
 \label{Phi_m_theta_star_partial_F_n_partial_2_Phi} 
  &\leq C \Big(  \Big \| \frac{1}{\sqrt{n}} \dfrac{\partial \Phi_n(\theta^*)}{\partial \theta} \Big \|^2
  +  \Big \| \frac{1}{\sqrt{n}}  \Big( \dfrac{\partial\widehat{\Phi}_n(\theta^*)}{\partial \theta}  - \dfrac{\partial \Phi_n(\theta^*)}{\partial \theta} \Big) \Big \|^2 \Big) \| \widehat{F}_n^{-1}(m) \|^2 \Big\| \dfrac{1}{n} \dfrac{\partial^2 \widehat{\Phi}_n\big( \overline{\theta}(m) \big)}{\partial \theta^2} \Big \|  .
\end{align}
According to the law of the law of iterated logarithm of martingales (see Stout (1970, 1974)), one can find a constant $C>0$ such that
\begin{equation} \label{Phi_m_theta_star_LIL} 
\limsup_{n \rightarrow \infty} \displaystyle{ \frac{1}{\sqrt{n \log \log n }} \Big\| \dfrac{\partial\Phi_n(\theta^*)}{\partial \theta} \Big\| = \limsup_{n \rightarrow \infty}  \dfrac{1}{\sqrt{n \log \log n}} \Big\| \sum_{t=1}^n \dfrac{\partial\varphi_t(\theta^*)}{\partial \theta} } \Big\| = C ~a.s.
\end{equation}
Hence, from (\ref{Phi_m_theta_star_partial_F_n_partial_2_Phi}), (\textbf{A3}) with (\ref{eq_assump_app_dif_log2}), (\ref{F_n_partial2_widehat_Phi_conv}) it holds that  
\begin{equation} \label{Phi_m_theta_star_O_1_loglog}
\limsup_{n \rightarrow \infty} \frac{1}{ \log \log n } \big| \widehat{\Phi}_{n}\big(\widehat{\theta}(m)\big)-\widehat{\Phi}_{n}\big(\theta^*) \big| \leq C ~ a.s.
\end{equation}  
for some constant $C>0$.
By using the same arguments with $m=m^*$, it follows that,
\begin{equation} \label{Phi_m_star_theta_star_O_1_loglog}
\limsup_{n \rightarrow \infty} \frac{1}{ \log \log n } | \widehat{\Phi}_{n}\big(\widehat{\theta}(m^*)\big)-\widehat{\Phi}_{n}\big(\theta^*) |  \leq C ~ a.s.
\end{equation}
for some $C>0$.
Thus, (\ref{Phi_m_star_Phi_m_limsup}) follows from (\ref{Phi_m_theta_star_O_1_loglog}) and (\ref{Phi_m_star_theta_star_O_1_loglog}).  
By virtue of (\ref{C_m_star_C_m_ps}) and (\ref{Phi_m_star_Phi_m_limsup}), we get $a.s$
\begin{align*}
 \nonumber \liminf_{n \rightarrow \infty} \dfrac{1}{\log \log n} \big( \widehat{C}(m) - \widehat{C}(m^*) \big) &\geq \liminf_{n \rightarrow \infty} \Big( -\frac{1}{ \log \log n } \big| \widehat{\Phi}_n\big(\widehat{\theta}(m)\big)-\widehat{\phi}_n\big(\widehat{\theta}(m^*)\big) \big| \Big) +   ( |m|-|m^*|)  \liminf_{n \rightarrow \infty} \dfrac{\kappa_n}{\log \log n} \\
 &\geq -C +  ( |m|-|m^*|)  \liminf_{n \rightarrow \infty} \dfrac{\kappa_n}{\log \log n}. 
\end{align*}
Therefore, since $m  \supsetneq  m^*$, if $\liminf_{n \rightarrow \infty} \dfrac{\kappa_n}{\log \log n} > 2C$, then 
\[ \liminf_{n \rightarrow \infty} \dfrac{1}{\log \log n} \big( \widehat{C}(m) - \widehat{C}(m^*) \big) > C > 0 ~ a.s. \]
Hence, 
\begin{equation} \label{C_m_C_m_star}
 \widehat{C}(m) - \widehat{C}(m^*) > 0 ~ a.s. ~ \text{  for } n \text{  large enough}. 
\end{equation}
From the proof of the first part of Theorem \ref{th1_mce}, one can see that (\ref{C_m_C_m_star}) also holds when $m  \nsupseteq  m^*$. Thus,
$ \widehat{m}= \underset{m \in \mathcal{M}}{\text{argmin}} ~ \widehat{C}(m) = \underset{m \in \mathcal{M}}{\text{argmin}} \big( \widehat{C}(m) - \widehat{C}(m^*) \big) \limitepsn m^*$; which completes the proof of the first part. 

\medskip

\noindent 2. Assume there exists $\widetilde{m} \in \mathcal{M}$ such that $m^* \subsetneq \widetilde{m} $. 
 In the proof of this second part, we set $\partial \psi (\theta) / \partial \theta := (\partial \psi(\theta) / \partial \theta_i)_{i \in \widetilde{m}}$ for any function $\psi$ and $\theta \in \Theta$.
 %
%
 From (\ref{Phi_m_theta_star_partial_F_n}) and (\textbf{A3}) with (\ref{eq_assump_app_dif_log2}) and (\ref{Phi_m_theta_star_LIL} ), we obtain
 \begin{align}
\nonumber &\dfrac{1}{\log \log n} \left( \widehat{\Phi}_n\big(\theta^*)  - \widehat{\Phi}_n\big(\widehat{\theta}(\widetilde{m})\big) \right) \\
\nonumber &=   \Big( \frac{1}{\sqrt{ 2n \log \log n}} \dfrac{\partial\widehat{\Phi}_n(\theta^*)}{\partial \theta^T} \Big)  \widehat{F}_n^{-1}(\widetilde{m}) \Big( \dfrac{1}{n} \dfrac{\partial^2 \widehat{\Phi}_n\big( \overline{\theta}(\widetilde{m}) \big)}{\partial \theta^2} \Big) \widehat{F}_n^{-1}(\widetilde{m}) \Big( \frac{1}{\sqrt{ 2n \log \log n}} \dfrac{\partial\widehat{\Phi}_n(\theta^*)}{\partial \theta} \Big) \\
\nonumber
&=   \Big( \frac{1}{\sqrt{2n \log \log n}} \dfrac{\partial \Phi_n(\theta^*)}{\partial \theta^T}  + o(1) \Big) \big( F(\widetilde{m})^{-1} + o(1) \big) \Big( \frac{1}{\sqrt{ 2n \log \log n }} \dfrac{\partial \Phi_n(\theta^*)}{\partial \theta}  + o(1) \Big) \\
\nonumber
&= \Big( \frac{1}{\sqrt{2n \log \log n}} \dfrac{\partial \Phi_n(\theta^*)}{\partial \theta^T}\Big) F(\widetilde{m})^{-1}  \Big( \frac{1}{\sqrt{ 2n \log \log n }} \dfrac{\partial \Phi_n(\theta^*)}{\partial \theta} \Big)  + o(1) \\
\nonumber
&= \dfrac{1}{\vartheta} \Big( \frac{1}{\sqrt{2n \log \log n}} \dfrac{\partial \Phi_n(\theta^*)}{\partial \theta^T}\Big) G(\widetilde{m})^{-1}  \Big( \frac{1}{\sqrt{ 2n \log \log n }} \dfrac{\partial \Phi_n(\theta^*)}{\partial \theta} \Big)  + o(1) \\
\nonumber
&=  \dfrac{1}{\vartheta} \left\| \frac{1}{\sqrt{2n \log \log n}}  G(\widetilde{m})^{-1/2} \dfrac{\partial \Phi_n(\theta^*)}{\partial \theta}  \right \|^2 + o(1) \\
\label{Phi_m_theta_star_partial_ps}
&=  \dfrac{1}{\vartheta}   \frac{1}{ 2n \log \log n} \sum_{i=1}^{|\widetilde{m}|} \zeta_{n,j}^2(\theta^*)  + o(1), 
\end{align}
where $\zeta_{n,j}(\theta^*)$ is a sum of $n$ term  of a stationary ergodic martingale difference sequence with unit variance.
 By virtue of the law of the iterated logarithm, it holds for any $j=1,\cdots,|\widetilde{m}|$
\begin{equation}\label{LIL_limsup_liminf}
 \limsup_{n \rightarrow \infty}  \frac{1}{\sqrt{2n \log \log n}}  \zeta_{n,j}(\theta^*)  = 1 ~  \text{ and } ~   \liminf_{n \rightarrow \infty}  \frac{1}{\sqrt{2n \log \log n}}  \zeta_{n,j}(\theta^*)  = -1 . 
 \end{equation}
Hence, (\ref{Phi_m_theta_star_partial_ps}) implies
\begin{equation}\label{Phi_theta_star_m_tilde_ps_log_log_vartheta}
 \dfrac{1}{\log \log n} \left( \widehat{\Phi}_n\big(\theta^*)  - \widehat{\Phi}_n\big(\widehat{\theta}(\widetilde{m})\big) \right) \limitepsn \dfrac{|\widetilde{m}|}{\vartheta}.
\end{equation}
 By going along similar lines as above and using (\ref{LIL_limsup_liminf}), we obtain
\begin{equation} \label{Phi_theta_star_m_star_ps_log_log_vartheta}
  \dfrac{1}{\log \log n} \left( \widehat{\Phi}_n\big(\theta^*)  - \widehat{\Phi}_n\big(\widehat{\theta}(m^*)\big) \right) \limitepsn \dfrac{|m^*|}{\vartheta}.
\end{equation} 
 From (\ref{Phi_theta_star_m_tilde_ps_log_log_vartheta}) and (\ref{Phi_theta_star_m_star_ps_log_log_vartheta}), we get
 \begin{align}
 \nonumber  \dfrac{1}{\log \log n} \left( \widehat{\Phi}_n\big(\widehat{\theta}(\widetilde{m})\big) - \widehat{\Phi}_n\big(\widehat{\theta}(m^*)\big) \right) &=    \dfrac{1}{\log \log n} \left( \widehat{\Phi}_n\big(\theta^*\big) - \widehat{\Phi}_n\big(\widehat{\theta}(m^*)\big)  \right) 
 +   \dfrac{1}{\log \log n} \left( \widehat{\Phi}_n\big(\widehat{\theta}(\widetilde{m})\big) - \widehat{\Phi}_n\big(\theta^*\big) \right) \\
 \label{Phi_m_tilde_m_star_ps_log_log_vartheta}
 & \limitepsn \dfrac{|m^*|}{\vartheta} - \dfrac{|\widetilde{m}|}{\vartheta} . 
\end{align}   
Thus, according to (\ref{C_m_star_C_m_ps}), we have
\begin{align*}
  \limsup_{n \rightarrow \infty} \dfrac{1}{\log \log n} \big( \widehat{C}(\widetilde{m}) - \widehat{C}(m^*)  \big) &\leq \limsup_{n \rightarrow \infty} \dfrac{1}{\log \log n} \left( \widehat{\Phi}_n\big(\widehat{\theta}(\widetilde{m})\big) - \widehat{\Phi}_n\big(\widehat{\theta}(m^*)\big) \right)  +  ( |\widetilde{m}|-|m^*|)  \limsup_{n \rightarrow \infty} \dfrac{\kappa_n }{\log \log n} \\
  &\leq  \dfrac{ |m^*| - |\widetilde{m}|}{\vartheta} + ( |\widetilde{m}|-|m^*|)  \limsup_{n \rightarrow \infty} \dfrac{\kappa_n }{\log \log n} 
 \end{align*}
Since $|m^*| < |\widetilde{m}| $, if $ \limsup_{n \rightarrow \infty}(\kappa_n / \log \log n) < 1/\vartheta$ then we can find a constant $C>0$, such that
\[   \limsup_{n \rightarrow \infty} \dfrac{1}{\log \log n} \big( \widehat{C}(\widetilde{m}) - \widehat{C}(m^*)  \big) \leq -C <0 ~ a.s.  \]
Consequently, $\widehat{C}(\widetilde{m}) - \widehat{C}(m^*) \limitepsn - \infty$; that is, $m^*$ is not asymptotically preferred, which shows that the strong consistency fails.     

\medskip

\noindent Now, assume that  $m^* \not\subset m$ for all $m \in \mathcal{M} \setminus \{m^*\}$.
From the proof of the first part of Theorem \ref{th1_mce}, one can see that (\ref{C_m_C_m_star}) holds when $m  \nsupseteq  m^*$.
That is, (\ref{C_m_C_m_star}) holds for all $m \in \mathcal{M} \setminus \{m^*\}$. Which shows that the strong consistency holds in this case, and completes the proof of the theorem.
\begin{flushright}
$\blacksquare$ 
\end{flushright}
%
 
\bibliographystyle{acm}

\end{document}